\documentclass[11pt,reqno]{amsart}

\usepackage[text={160mm,240mm},centering]{geometry}            
\usepackage{amssymb,amsmath,setspace,geometry,indentfirst,changepage,inputenc,amsthm}
\usepackage{mathrsfs,amsfonts,savesym,graphicx,bm,color}

\usepackage{graphicx}
\usepackage{amssymb}
\usepackage{dsfont}

\usepackage{lscape}

\usepackage{arydshln}

\usepackage{tikz}
\usetikzlibrary{positioning}

\usepackage{float}

\usepackage[all,2cell]{xy}    
\UseAllTwocells               

\usepackage{amsmath,amsthm}

\usepackage[mathscr]{eucal}
\usepackage[all]{xy}
\usepackage{mathrsfs}
\usepackage{hyperref}
\usepackage{color}

\newtheorem{theorem}{Theorem}[section]
\newtheorem{lemma}[theorem]{Lemma}

\newtheorem{conjecture}[theorem]{Conjecture}
\newtheorem{corollary}[theorem]{Corollary}
\newtheorem{proposition}[theorem]{Proposition}

\newtheorem{definition-lemma}[theorem]{Definition-Lemma}
\newtheorem{definition-theorem}[theorem]{Definition-Theorem}

\newtheorem*{conj}{Conjecture}

\theoremstyle{definition}
\newtheorem{example}[theorem]{Example}
\newtheorem{definition}[theorem]{Definition}

\newtheorem{remark}[theorem]{Remark}

\newtheorem*{ack}{Acknowledgements}

\usepackage{fancyhdr}
\allowdisplaybreaks

\title{On Stringy E-functions and the Non-negativity Conjecture for Determinantal Varieties}

\author{Yifan Chen}
\address{
	Department of Mathematical Sciences,
	Tsinghua University,
	Beijing, 100084, P. R. China.}
\email{c-yf20@tsinghua.org.cn}

\author{Huaiqing Zuo}
\address{Department of Mathematical Sciences,
	Tsinghua University,
	Beijing, 100084, P. R. China.}
\email{hqzuo@mail.tsinghua.edu.cn}

\begin{document}
	
	\maketitle
	
	%
	
	\begin{abstract}
We compute the stringy E-functions of determinantal varieties and establish that the stringy E-function of a determinantal variety coincides with the E-function of the product of a Grassmannian and an affine space. Furthermore, a similar result holds for the projectivization of determinantal varieties. As an application, we verify the non-negativity conjecture for the stringy Hodge numbers of these varieties.
	\end{abstract}
	
	\tableofcontents
	
	\section{Introduction}
In \cite{Batyrev_Stringy_Invariant} and \cite{Stringy_Hodge_Number}, Batyrev introduced the "stringy E-function" as an invariant of singularities. It is a rational function in variables $u$ and $v$, and in general, it is only defined for log terminal and $\mathbb{Q}$-Gorenstein singularities. Batyrev utilized this function to generalize the E-polynomial and study mirror symmetry for Calabi-Yau varieties. He also proved a version of the McKay correspondence using the stringy E-function. Although the stringy E-function can be expressed in terms of resolution data,  its calculation  is generally challenging, and explicit expressions are known only in a few special cases.  For instance, Roczen and Dais derived  the formula for the stringy E-function of 3-dimensional ADE singularities in \cite{Stringy_Invariant_of_ADE_Singularities}, and the latter also computed it for a class of absolutely isolated singularities in \cite{Stringy_Invariant_of_Absolutely_Isolated_Singularities}. Schepers provided the explicit form for strictly canonical non-degenerate singularities  in \cite{Stringy_Invariant_of_Strictly_Canonical_Nondegenerate_Singularities}  and, together with Veys,   studied the case of Brieskorn singularities in \cite{Stringy_E-function_of_Hypersurface}.    

Determinantal varieties are important examples of varieties to study. In this paper, we compute the stringy E-functions of these varieties.  First,  we fix the notation. Let $\mathcal{M}:=\mathbb{A}^{rs}$ denote the space of $r \times s$ matrices, and let $D^{k}$ be the variety of matrices with rank $\le k$,  which  is defined by the vanishing of $k+1$ minors.  As shown in Section $8$ of  \cite{Determinantal_Ring},  $D^{k}$ is Gorenstein if and only if $r=s$.  Therefore,  we only  consider the case $r=s$.  We now state one of our main results as follows.

\begin{theorem}(Theorem \ref{formula})\label{th1}
\textit{With the notation as above, $D^{k}$ has Gorenstein canonical singularities, and its stringy E-function is given by}
\begin{flalign*}
	E_{st}(D^{k})=E(\mathbb{L}^{kr}\prod_{i=1}^{k}\frac{\mathbb{L}^{i+r-k}-1}{\mathbb{L}^{i}-1})=E(\mathbb{L}^{kr} \cdot G(k,r)),
\end{flalign*}
\textit{where $G(k,r)$ denotes the Grassmannian of $k$-dimensional subspaces in an $r$-dimensional vector space, $E_{st}(\cdot)$ represents  the stringy E-function,  and $E(\cdot)$ denotes the E-polynomial(see Section \ref{E-polynomial_and_stringy_E-function} for the definition in details).}
\end{theorem}

In \cite{Stringy_Hodge_Number}, Batyrev also defined "stringy Hodge number" using stringy E-function, generalizing the usual Hodge number. The definition is as follows.
\begin{definition}
Suppose $X$ is a $\mathbb{Q}$-Gorenstein variety with at worst log terminal singularity. If the stringy E-function of $X$ is a polynomial, written as $$E_{st}(X)=\sum_{p,q}(-1)^{p+q}h^{p,q}u^{p}v^{q},$$ then we define these $h^{p,q}$ to be the stringy Hodge numbers of $X$. If the stringy E-function of $X$ is not a polynomial, then we say the stringy Hodge number of $X$ does not exist.
\end{definition}

The classical properties of Hodge numbers for smooth projective varieties
are still satisfied by Batyrev’s stringy Hodge numbers for Gorenstein canonical
projective varieties with a polynomial stringy E-function. But a serious problem appears: it is not clear that they are nonnegative.

\begin{conj}(\cite{Stringy_Hodge_Number},  Conjecture 3.10). Let $X$ be a Gorenstein canonical projective variety. Assume that $E_{s t}(X)$ is a polynomial. Then all stringy Hodge numbers $h_{s t}^{p, q}(X)$ are nonnegative.
\end{conj}


This is the famous Batyrev’s "non-negativity" conjecture and it is one of the most important problems in this field. In practice, people guess that these stringy Hodge numbers are actually the dimensions of some cohomology groups, so they are always nonnegative. There are only some partial results about this conjecture, for example, Yasuda proved the quotient singularity case in \cite{Non-negetivity_Conjecture_for_Quotient_singularity}, and in \cite{Non-negativity_Conjecture_for_Toric_Artin_Stack} the case of some special Artin toric stack is proven by Satriano and others. In this paper, we verify  the  "non-negativity" conjecture for the projectivization of determinantal varieties by computing the stringy E-function of it. To be precise, let $\hat{D}^{k}$ be the projectivization of $D^{k}$ in $\mathcal{M}$, then we obtain the following result. 

\begin{theorem}\label{th2}
$\hat{D}^{k}$ has Gorenstein canonical singularities,  and its stringy E-function is given as follows
\begin{flalign*}
	E_{st}(\hat{D}^{k}_{r,r})=E(\frac{\mathbb{L}^{kr}-1}{\mathbb{L}-1} \cdot \prod_{i=1}^{k}\frac{\mathbb{L}^{i+r-k}-1}{\mathbb{L}^{i}-1})=\frac{(uv)^{kr}-1}{uv-1} \cdot E(G(k,r)).
\end{flalign*}
In particular, the Batyrev’s "non-negativity" conjecture holds for $\hat{D}^{k}$.
\end{theorem}

The paper is organized as follows. Section 2 provides preliminary materials.  In Section 3, we investigate the singularities of determinantal varieties. Section 4 and 5 present the proofs  of Theorem \ref{th1} and Theorem \ref{th2},  respectively.

\begin{ack}
Y. Chen and H. Zuo were supported by BJNFS Grant 1252009.  H. Zuo was supported by NSFC Grant 12271280.
\end{ack}

	\section{Preliminary}
	\subsection{Jet Scheme, Grothendieck Ring and Motivic Integration}\label{sec_Grothendieck_Ring}
In this section, we outline the foundational concepts of jet schemes, the Grothendieck ring, and motivic integration.
\begin{definition}[Jet Scheme]
	Let $X$ be a scheme over a field $k$, for every $m \in \mathbb{N}$, consider the functor from $k$-schemes to set
	\begin{flalign*}
		Z \mapsto \mathrm{Hom}(Z \times_{\mathrm{Spec}(k)} \mathrm{Spec} (k[t]/(t^{m+1})),X).
	\end{flalign*}
	There is a $k$-scheme represents this functor, which is called the $m$-th jet scheme of $X$, denoted by $J_{m}(X)$(see Theorem 2.1, \cite{Jet_Schemes_Ishii}), i.e.
	\begin{flalign*}
		\mathrm{Hom}(Z,J_{m}(X))=\mathrm{Hom}(Z \times_{\mathrm{Spec}(k)} \mathrm{Spec} (k[t]/(t^{m+1})),X).
	\end{flalign*}
\end{definition}

For $1 \le i \le j$, the truncation map $k[t]/(t^{j}) \to k[t]/(t^{i})$ induces natural projections between jet schemes $\psi_{i,j}: J_{j}(X) \to J_{i}(X)$. If we identify $J_{0}(X)$ with $X$, we can get natural projection $\pi_{m}: J_{m}(X) \to X$. One can check that $\{J_{m}(X)\}_{m}$ is an inverse system and the inverse limit $J_{\infty}(X):= \varprojlim_m J_{m}(X)$ is called the arc space of $X$, inducing natural projections $\psi_{i}: J_{\infty}(X) \to J_{i}(X)$. If $X$ is of finite type, $J_{m}(X)$ is also of finite type for each $m \in \mathbb{N}$, but $J_{\infty}(X)$ is usually not.

Next we introduce the Grothendieck ring, which makes the set of all algebraic varieties an abelian group.
\begin{definition}[Grothendieck Ring]\label{Ch2_Sec2_Grothendieck_Ring}
	Let $k$ be a field and let $\mathrm{Var}_{k}$ be the category of $k$-varieties. The Grothendieck group of $k$-varieties, $K_0[\mathrm{Var}_k]$, is defined to be the quotient group of the free abelian group with basis $\{[X]\}_{X\in\mathrm{Var}_k}$, modulo the following relations.
	\begin{align*}
		& [X]-[Y], \ X \simeq Y\\
		& [X]-[X_{\mathrm{red}}]\\
		& [X]-[U]-[X\setminus U],\ U\subseteq X\ \mathrm{open}.
	\end{align*}
	One can further define a multiplication structure on $K_0[\mathrm{Var}_k]$ by
	\begin{displaymath}
		[X] \cdot [Y] := [X\times Y].
	\end{displaymath}
	This makes $K_0[\mathrm{Var}_k]$ a ring, called the Grothendieck ring of $k$-varieties.
	
	Let $\mathbb L = [\mathbb A_k^1]$ and $K_0[\mathrm{Var}_k]_{\mathbb L}$ be the localization at $\mathbb L$.
\end{definition}

\noindent\textbf{Kontsevich's Completion}

Kontsevich's completion allows us to take limit in Grothendieck ring, and we are going to introduce this. There is a natural decreasing filtration $F^\bullet$ on $K_0[\mathrm{Var}_{k}]_{\mathbb L}$. For $m\in \mathbb Z$, $F^m$ is the subgroup of $K_0[\mathrm{Var}_{k}]_{\mathbb L}$ generated by $[S] \cdot \mathbb L^{-i}$ with $\dim S- i \leq -m$. It is actually a ring filtration i.e. $F^m \cdot F^n \subseteq F^{m+n}$.
\begin{definition}
	The Kontsevich's completed Grothendieck ring $\widehat{K_0[\mathrm{Var}_{k}]}$ is defined to be 
	\begin{displaymath}
		\widehat{K_0[\mathrm{Var}_{k}]} := \varprojlim_{m\in \mathbb Z} K_0[\mathrm{Var}_{k}]_{\mathbb L}/F^m.
	\end{displaymath}
\end{definition}

Now we fix an algebraic variety $X$ of pure dimension $d$ and we will do the motivic integration on the space $J_{\infty}(X)$. Firstly we need to define the measurable sets, which turns out to be the cylinders in $J_{\infty}(X)$, with the following definition.
\begin{definition}[Cylinder]
	A subset $C \subset J_{\infty}(X)$ is called a cylinder if $C=\psi_{m}^{-1}(A)$ for some constructible subset $A \subset J_{m}(X)$. 
\end{definition}

We want to define the measure $\mu_{J_{\infty}(X)}$ on $J_{\infty}(X)$ such that $\mu_{J_{\infty}(X)}(C)=\frac{[\psi_{m}(C)]}{\mathbb{L}^{md}}$ for a cylinder $C=\psi_{m}^{-1}(A)$. To show that this definition is independent of the choice of $m$, we use the following theorem.

\begin{theorem}[\cite{Motivic_Integration_on_Arbitrary_Varieties_Denef_Loeser}]
	If $C \subset J_{\infty}(X)$ is a cylinder and $C \cap J_{\infty}(X_{\mathrm{sing}})=\emptyset$, then there exists integer $m$ such that
	
	(1) $\psi_{m}(C)$ is constructible and $C=\psi_{m}^{-1}(\psi_{m}(C))$.
	
	(2) For $n \ge m$, the projection $\psi_{n+1}(C) \to \psi_{n}(C)$ is a piecewise trivial fibration with fiber $\mathbb{A}^{d}$.
	
	In particular, $\frac{[\psi_{n}(C)]}{\mathbb{L}^{nd}}$ are equal for $n \ge m$.
\end{theorem} 

For general cylinder $C$, the condition of the above theorem may not hold and we will define the measure with value in Kontsevich's completed Grothendieck ring $\widehat{K_0[\mathrm{Var}_{k}]}$, with the following theorem.

\begin{theorem}[\cite{Motivic_Integration_on_Arbitrary_Varieties_Denef_Loeser}]
	Let $C$ be a cylinder of $J_{\infty}(X)$, then the limit 
	\begin{flalign*}
		\mu_{J_{\infty}(X)}(C):=\lim_{m \rightarrow \infty} \frac{[\psi_{m}(C)]}{\mathbb{L}^{md}}
	\end{flalign*}
	exists in $\widehat{K_0[\mathrm{Var}_{k}]}$.
\end{theorem}

Now we can define the motivic integration.
\begin{definition}
	Suppose $C \subset J_{\infty}(X)$ is a cylinder and $\alpha: C \to \mathbb{Z} \cup \{\infty\}$ with cylinder fiber, we define the motivic integration
	\begin{flalign*}
		\int_{C}\mathbb{L}^{-\alpha}d\mu_{J_{\infty}(X)}=\sum_{n \in \mathbb{Z}}\mu_{J_{\infty}(X)}(\alpha^{-1}(n))\mathbb{L}^{-n},
	\end{flalign*}
	if the right hand side convergences.
\end{definition}

One of the most famous problem related to motivic integration is monodromy conjecture, which gives the relation between the poles of motivic zeta function and the roots of Bernstein-Sato polynomial of a polynomial over complex field. The definition of the motivic zeta function is given by the motivic integration as follows.

\begin{definition}
Let $f \in \mathbb{C}[x_{1},...,x_{n}]$ be a polynomial over complex field, then the motivic zeta function of $f$ is given by
\begin{flalign*}
	Z_{f}^{mot}(s):=\int_{J_{\infty}(\mathbb{C}^{n})}\mathbb{L}^{-\mathrm{ord}_{t}(f)}d\mu_{J_{\infty}(\mathbb{C}^{n})},
\end{flalign*}
where $\mathrm{ord}_{t}(f)$ is a function sending an arc $\psi \in J_{\infty}(\mathbb{C}^{n})$ to the order of the formal power series $f(\psi)$.
\end{definition}

For the introduction of Bernstein-Sato polynomial, one can refer to, for example, \cite{Bernstein_Sato_Notes}.

\begin{conjecture}[Monodromy conjecture]\label{monodromy_conjecture}
Let $f \in \mathbb{C}[x_{1},...,x_{n}]$ be a non-constant polynomial over complex field, then the poles of $Z_{f}^{mot}(s)$ are the roots of $b_{f}(s)$, the Bernstein-Sato polynomial of $f$.
\end{conjecture}

In \cite{Motivic_Zeta_Function_on_Q_Gorenstein_Varieties_Leon_Martin_Veys_Viu_Sos}, the authors gave another measure $\mu^{\mathbb{Q}-Gor}$ on $J_{\infty}(X)$ for those $X$ which are $\mathbb{Q}$-Gorenstein. They used this measure to give a formula of motivic zeta function via embedded $\mathbb{Q}$-resolution, which was utilized to prove the monodromy conjecture of semi-quasihomogeneous polynomials by  Blanco,  Budur and  van der Veer in \cite{Monodromy_Conjecture_for_Semi-Quasihomogeneous_Hypersurfaces_Blanco_Budur_vdV}. In this paper, we will use this measure to study stringy E-functions. To establish this, some preliminary groundwork is required.

\begin{definition}[Definition $1.1$ in \cite{Motivic_Zeta_Function_on_Q_Gorenstein_Varieties_Leon_Martin_Veys_Viu_Sos}]
Let $X$ be a complex variety of pure dimension $d$ and $\Omega_{X}^{d}$ be the sheaf of $d$-th exterior product of differential form on $X$. We fix an arc $\psi \in J_{\infty}(X)-J_{\infty}(X_{sing})$ and we regard it as a morphism $\psi: \mathrm{Spec}\mathbb{C}[[t]] \longrightarrow X$. Now we consider $\psi^{*}(\Omega_{X}^{d})$ as a $\mathbb{C}[[t]]$-module and $V:=\psi^{*}(\Omega_{X}^{d})\otimes_{\mathbb{C}[[t]]}\mathbb{C}((t))$ as a $\mathbb{C}((t))$-vector space. We denote $L_{X}$ as the image of $\psi^{*}(\Omega_{X}^{d})$ in $V$. 

Suppose $\omega$ is an invertible $\mathcal{O}_{X}$-subsheaf of $\Omega_{X}^{d} \otimes_{\mathbb{C}}\mathbb{C}(X)$ and we denote $\Lambda_{X}$ as the image of $\psi^{*}(\omega)$ in $V$. If $\Lambda_{X}=0$, then we define $\mathrm{ord}_{t}\omega(\psi)=\infty$. Otherwise we have $\Lambda_{X}=t^{e}L_{X}$ for some integer $e$, and we define $\mathrm{ord}_{t}\omega(\psi)=e$. 

More generally, if $\omega$ is an invertible $\mathcal{O}_{X}$-subsheaf of $(\Omega_{X}^{d})^{\otimes r} \otimes_{\mathbb{C}}\mathbb{C}(X)$, we denote $\Lambda_{X}$ as the image of $\psi^{*}(\omega)$ in $V^{\otimes r}$ and we have $\Lambda_{X}=t^{e}L_{X}^{\otimes r}$ for some $e \in \mathbb{Z} \cup \{\infty\}$ in $V^{\otimes r}$. We set $\mathrm{ord}_{t}\omega(\psi)=e$.
\end{definition}

Now we suppose $X$ is a $\mathbb{Q}$-Gorenstein variety over $\mathbb{C}$ of pure dimension $d$ with at worst log terminal singularity. We define $\omega_{X}:=i_{*}(\Omega_{X_{reg}}^{d})$, where $i: X_{reg} \longrightarrow X$ is the inclusion of the smooth part of $X$ into $X$. Since $X$ is $\mathbb{Q}$-Gorenstein, there exists an integer $r$ such that $\omega_{X}^{[r]}:=i_{*}((\Omega_{X_{reg}}^{d})^{\otimes r})$ is an invertible sheaf. 

\begin{definition}[Motivic $\mathbb{Q}$-Gorenstein measure, Definition $1.4$ in \cite{Motivic_Zeta_Function_on_Q_Gorenstein_Varieties_Leon_Martin_Veys_Viu_Sos}]
For a cylinder $C \subset J_{\infty}(X)$, we define the motivic $\mathbb{Q}$-Gorenstein measure of it to be
\begin{flalign*}
\mu^{\mathbb{Q}-Gor}_{J_{\infty}(X)}(C):=\int_{C}\mathbb{L}^{-\frac{1}{r} \mathrm{ord}_{t}\omega_{X}^{[r]}}d\mu_{J_{\infty}(X)}.
\end{flalign*}
\end{definition}

\begin{remark}\label{two_measure_coincide}
If $X$ is smooth, $i: X_{reg} \longrightarrow X$ is identity and $\omega_{X}=\Omega_{X_{reg}}^{d}$, hence we have $\mu^{\mathbb{Q}-Gor}_{J_{\infty}(X)}=\mu_{J_{\infty}(X)}$, i.e. the usual measure and $\mathbb{Q}$-Gorenstein measure coincide when $X$ is smooth. 
\end{remark}

\subsection{E-polynomials and stringy E-functions}\label{E-polynomial_and_stringy_E-function}
In this subsection we introduce the definition of E-polynomial and stringy E-function. Suppose $X$ is a complex $\mathbb{Q}$-Gorenstein variety with at worst log terminal singularities of dimension $n$. Deligne proved in \cite{Hodge_Theory2}, \cite{Hodge_Theory} that there is a natural mixed Hodge structure of cohomology group $H^{k}(X, \mathbb{Q})$. This is also true for cohomology group with compact support $H^{k}_{c}(X,\mathbb{Q})$. That is, an increasing weighted filtration $$0=W_{-1} \subset W_{0} \subset ... \subset W_{2k}=H^{k}_{c}(X,\mathbb{Q})$$ and a decreasing Hodge filtration $$0=F^{k+1} \subset ... \subset F^{0}=H^{k}_{c}(X,\mathbb{C})$$ such that there is a pure Hodge structure of weight $l$ on the graded quotient $\mathrm{Gr}_{l}^{W}H^{k}_{c}(X):=W_{l}/W_{l-1}$ induced by the filtration $F^{\bullet}$. Let $F^{p}\mathrm{Gr}_{l}^{W}H^{k}_{c}(X)$ be the complexified image of $F^{p} \cap W_{l}$ in $W_{l}/W_{l-1} \otimes \mathbb{C}$, then we have Hodge-Deligne numbers $$h^{p,q}(H^{k}_{c}(X,\mathbb{C})):=\mathrm{dim}_{\mathbb{C}}(F^{p}\mathrm{Gr}_{p+q}^{W}H^{k}_{c}(X) \cap \overline{F^{q}\mathrm{Gr}_{p+q}^{W}H^{k}_{c}(X)}).$$
\begin{definition}[E-polynomial]
	The E-polynomial of $X$ is defined to be
	\begin{flalign*}
		E(X):=\sum_{0 \le p,q \le n}\sum_{0 \le k \le 2n}(-1)^{k}h^{p,q}(H^{k}_{c}(X,\mathbb{C}))u^{p}v^{q} \in \mathbb{Z}[u,v].
	\end{flalign*}
	In particular, if we let $u=v=1$, we can get the topological Euler number $e(X)$.
\end{definition}
Now we can introduce the definition of stringy E-function.
\begin{definition}[Stringy E-function]\label{def_of_stringy_E-function}
	Suppose $X$ is a $\mathbb{Q}$-Gorenstein and log terminal variety of dimension $n$. Take a log resolution $\phi: Y \longrightarrow X$, and let $D_{i}(i \in S)$ be its irreducible components of relative canonical divisor and $a_{i}(i \in S)$ be its log discrepancies. For $I \subset S$, we denote $\mathring{D_{I}}:=\bigcap_{i \in I}D_{i} \backslash \bigcup_{i \notin I}D_{i}$, then the stringy E-function, i.e. the stringy E-function of $X$ is defined to be
	\begin{flalign*}
		E_{st}(X):=\sum_{I \subset S}E(\mathring{D_{I}})\prod_{i \in I}\frac{uv-1}{(uv)^{a_{i}}-1},
	\end{flalign*}
	where $E(\cdot)$ means the E-polynomial.
	Using Theorem $2.15$ in \cite{Introduction_to_Motivic_Integration}, we can rewrite the stringy E-function of $X$ into the following integration
	\begin{flalign*}
		E_{st}(X)=E(\int_{J_{\infty}(Y)}\mathbb{L}^{-\mathrm{ord}_{t}K_{\phi}}d\mu_{J_{\infty}(Y)}),
	\end{flalign*}
	where $K_{\phi}=K_{Y}-\phi^{*}(K_{X})$ is the relative canonical divisor.
\end{definition}

\subsection{Determinantal varieties}
In this subsection, we establish notations and review fundamental results concerning determinantal varieties and their arc spaces, primarily adopting the conventions for determinantal varieties as outlined  in \cite{Arcs_on_Determinantal_Variteties_Roi}. More results concerning the singularities of determinantal varieties can be found in, for example, Section $4$ in \cite{Local_Structure_of_Theta_Divisors} and Section $7$ in \cite{Deformations_with_cohomology_constraints}.

Let $\mathcal{M}=\mathbb{A}^{rs}$$(r \le s)$ be the space of $r \times s$ matrices and $D^{k}$$(0 \le k \le r-1)$ be the scheme of matrices with rank $\le k$, so $D^{k}$ is defined by the ideal of $(k+1)$-minors. Let $x_{i,j}(1 \le i \le r, 1 \le j \le s)$ be the elements of the matrix, so the coordinate ring of $\mathcal{M}$ is $k[x_{i,j}|1 \le i \le r, 1 \le j \le s]$. Assume $G=\mathrm{GL}_{r} \times \mathrm{GL}_{s}$ and we define the action of $G$ on $\mathcal{M}$ as follows.
\begin{flalign*}
	\mathrm{GL}_{r} \times \mathrm{GL}_{s} \times \mathcal{M} &\longrightarrow \mathcal{M}   \\
	(g,h) \cdot A &\mapsto gAh^{-1}.
\end{flalign*}
This action induces the action of $J_{\infty}(G)$ on $J_{\infty}(\mathcal{M})$ by just replacing the elements in matrices by the formal power series. Now $J_{\infty}(\mathcal{M})$ is decomposed into some orbits by this action, and every orbit has a standard form, as the following theorem in \cite{Arcs_on_Determinantal_Variteties_Roi}.
\begin{theorem}[Proposition 3.2 in \cite{Arcs_on_Determinantal_Variteties_Roi}]\label{Orbit decomposition}
Suppose $\infty \ge \lambda_{1} \ge \lambda_{2} \ge ... \ge \lambda_{r} \ge 0$ with $\lambda_{1},...,\lambda_{r} \in \mathbb{N}$, and $\lambda=(\lambda_{1},...,\lambda_{r})$. Then every standard form of the orbit is of the form $\delta_{\lambda}$, defined as follows:
\begin{flalign*}
	\delta_{\lambda}=\begin{pmatrix}
		t^{\lambda_{1}} & 0 & \dots & 0 & \dots & 0 \\
		0 & t^{\lambda_{2}} & \dots & 0 & \dots & 0 \\
		\vdots & \vdots & \ddots & \vdots & \dots & 0 \\
		0 & 0 & \dots & t^{\lambda_{r}} & \dots & 0
	\end{pmatrix}.
\end{flalign*}
Here we set $t^{\infty}=0$. So every orbit is indexed by some $\lambda$ and we denote the orbit corresponding to $\lambda$ by $\mathcal{C}_{\lambda}$. Moreover, $J_{\infty}(D^{k})$ is stable under the $J_{\infty}(G)$-action, and we have $\mathcal{C}_{\lambda} \subset J_{\infty}(D^{k})$ if and only if $\lambda_{1}=...=\lambda_{r-k}=\infty$.
\end{theorem}

\section{Singularities of determinantal varieties}
To define the stringy E-functions of determinantal varieties, we will verify in this section  whether $D^{k}$ satisfies the $\mathbb{Q}$-Gorenstein and log terminal conditions. Section $8$ in \cite{Determinantal_Ring} shows that $D^{k}$ is $\mathbb{Q}$-Gorenstein if and only if $r=s$, and actually when $r=s$, $D^{k}$ is Gorenstein, so we restrict ourselves to the case $r=s$. 

To verify the log terminal condition, it is necessary to analyze the resolution and the canonical divisor of $D^{k}$. Proposition $4.3$ in \cite{Minimal_Log_Discrepancy_of_Determinantal_Varieties} provides an explicit generator for the canonical divisor  of $D^{k}$, and we will introduce some notation to describe it. For $I,J \subset \{1,...,r\}$ with $|I|=|J|=l$, let $\Delta_{IJ}$ denotes the $l$-minor of $r \times r$ matrix formed by the rows and columns indexed by $I$ and $J$ respectively$(1 \le l \le r)$. Let $D(\Delta_{IJ})$ be the $\Delta_{IJ}$ nonzero part in $D^{k}$, so it is an open subscheme of $D^{k}$. It is well-known that the singular locus of $D^{k}$ is exactly $D^{k-1}$, so we have $D^{k}_{reg}=\bigcup_{|I|=|J|=k}D(\Delta_{IJ})$, where $D^{k}_{reg}$ denotes the smooth part of $D^{k}$. Let $\mathcal{S}_{IJ}=\{x_{i,j}\vert i \in I \ \mathrm{or} \ j \in J\}$. We define a differential $k(2r-k)$-form $\omega$ on $D(\Delta_{\{1,...,k\vert 1,...,k\}})$,
\begin{flalign*}
	\omega=\frac{1}{(\Delta_{\{1,...,k\vert 1,...,k\}})^{r-k}}\bigwedge_{x_{i,j} \in \mathcal{S}_{\{1,...,k\vert 1,...,k\}}}dx_{i,j}.
\end{flalign*}
Note that the dimension of $D^{k}$ is $k(2r-k)$, so $\omega$ is the differential form of the highest order. The next proposition shows that $\omega$ can be extended to the whole $D^{k}_{reg}$ and it is the generator of the sheaf $\omega_{D^{k}}=i_{*}(\Omega_{D^{k}_{reg}}^{k(2r-k)})$, where $i:D^{k}_{reg} \longrightarrow D^{k}$ is the inclusion.
\begin{proposition}[Proposition $4.3$ in \cite{Minimal_Log_Discrepancy_of_Determinantal_Varieties}]\label{Generator of canonical divisor}
	$\omega$ can be extended to a global section of $\omega_{D^{k}}$, and for each $D(\Delta_{IJ})$ we have
	\begin{flalign*}
		\omega|_{D(\Delta_{IJ})}=\pm \frac{1}{(\Delta_{IJ})^{r-k}}\bigwedge_{x_{i,j} \in \mathcal{S}_{IJ}}dx_{i,j}.
	\end{flalign*}
	Moreover, $\omega$ generates $\omega_{D^{k}}$.	
\end{proposition}

The embedded resolution of $D^{k}$ is given in \cite{Multiplier_Ideals_of_Determinantal_Ideal}, with the following proposition.
\begin{proposition}[Theorem $4.4$ in \cite{Multiplier_Ideals_of_Determinantal_Ideal}]\label{embedded resolution}
Let $\pi_{0}: A_{0} \longrightarrow A_{-1}:=\mathbb{C}^{r^{2}}$ be the blowup of $D^{0}$. Suppose $\pi_{0}:A_{0} \longrightarrow A_{-1},...,\pi_{i-1}:A_{i-1} \longrightarrow A_{i-2}$ have been defined, set $\pi_{i}: A_{i} \longrightarrow A_{i-1}$ to be the blowup of $\tilde{D^{i}}$ in $A_{i-1}$, where $\tilde{D^{i}}$ is the strict transform of $D^{i}$ via $\pi_{0} \circ \pi_{1} \circ ... \circ \pi_{i-1}:A_{i-1} \longrightarrow A_{-1}$. Then $\pi: A_{k} \longrightarrow A_{k-1} \longrightarrow ... \longrightarrow A_{0} \longrightarrow \mathbb{C}^{r^{2}}$ is an embedded resolution of $D^{k}$.
\end{proposition}
To give the resolution and log discrepancy of $D^{k}$, we also need the explicit process of these blowups. Let $x_{i,j}$ be the coordinates of $\mathbb{C}^{r^{2}}$, and $A_{0}=\{([y_{i,j}],(x_{i,j})) \in \mathbb{P}^{r^{2}-1} \times \mathbb{C}^{r^{2}}\; \vert\; \exists \,\theta, x_{i,j}=\theta y_{i,j}\}$ be the blowup of $0$. The exceptional divisor $E$ is defined by $\theta=0$, and the strict transform $\tilde{D^{k}}$ is given by the $(k+1)$-minors of $y_{i,j}$. $A_{0}$ is covered by $r^{2}$ charts, where each chart $U_{i,j}$ is given by $y_{i,j}=1$ for some $i,j$. Without loss of generality, we focus on the chart $U_{1,1}$. On $U_{1,1}$, we have $\theta,y_{i,j}((i,j) \neq (1,1))$ as new affine coordinates and $x_{1,1}=\theta, x_{i,j}=\theta y_{i,j}$ for $(i,j) \neq (1,1)$. In these coordinates, $\tilde{D^{k}}$ in $U_{1,1}$ is defined by the $(k+1)$-minors of the following matrix
\begin{flalign*}
	\begin{pmatrix}
		1       & y_{1,2} & \dots  & y_{1,r} \\
		y_{2,1} & y_{2,2} & \dots  & y_{2,r} \\
		\vdots  & \vdots  & \ddots & \vdots \\
		y_{r,1} & y_{r,2} & \dots  & y_{r,r}
	\end{pmatrix}.
\end{flalign*}
Let $f_{i,j}=y_{i,j}-y_{i,1}y_{1,j}$, where $2 \le i,j \le r$, then Lemma $3.13$ in \cite{Multiplier_Ideals_of_Determinantal_Ideal} shows that the ideal generated by $(k+1)$-minors of the matrix above is the same as the ideal generated by $k$-minors of the matrix $(f_{i,j})_{(r-1) \times (r-1)}$. Now we use $\theta,y_{1,2},...,y_{1,r},y_{2,1},...,y_{r,1},f_{i,j}(2 \le i,j \le r)$ as new coordinates, then $\tilde{D^{k}} \cap U_{1,1}$ is isomorphic to the product of $\mathbb{C}^{2r-1}$ and the determinantal variety with rank $\le k-1$ in $(r-1) \times (r-1)$ matrix. If we replace $k$ by $1$, we can get $\tilde{D^{1}}$ on each chart, which is exactly the center of the second blowup. This shows how these blowups in Proposition \ref{embedded resolution} work and that $\tilde{D^{i}}$ in $A_{i}$ are smooth. Now we can give the resolution of $D^{k}$ by just restricting our embedded resolution to our strict transform.
\begin{proposition}\label{resolution of Dk}
Assume the strict transform of $D^{k}$ in $A_{k-1}$ is $\tilde{D^{k}}$, then $\tilde{\pi}:=\pi_{0} \circ \pi_{1} \circ ... \circ \pi_{k-1}|_{\tilde{D^{k}}}: \tilde{D^{k}} \longrightarrow D^{k}$ is a resolution of $D^{k}$.
\end{proposition}
\begin{proof}
It is well-known that the singular locus of $D^{k}$ is $D^{k-1}$. By Proposition \ref{embedded resolution}, $\pi_{0} \circ \pi_{1} \circ ... \circ \pi_{k-1}$ is an embedded resolution of $D^{k-1}$, so its exceptional divisors $E_{0},...,E_{k-1}$ are simple normal crossings divisors. Also $\pi_{0} \circ \pi_{1} \circ ... \circ \pi_{k-1}$ is an isomorphism outside $D^{k-1}$, combining with the fact that $\tilde{D^{k}}$ is smooth we know that the assertion holds.
\end{proof}

We now compute the log discrepancy of $D^{k}$ using the resolution $\tilde{\pi}$.  Proposition \ref{Pr3.4} ensures that $D^{k}$  satisfies the necessary condition for defining  the stringy E-function.
\begin{proposition}\label{Pr3.4}
The log discrepancy of $D^{k}$ is $(k-i)(r-i)$, $i=0,1,...,k-1$. In particular, $D^{k}$ has canonical singularities. 
\end{proposition}
\begin{proof}
We denote the exceptional divisor of $\pi_{i}$ by $E_{i}$ and its intersection with the strict transform of $D^{k}$ by $\tilde{E_{i}}$. Firstly we compute the pull-back of $\omega$ via $\pi_{0}$, where $\omega$ is defined in Proposition \ref{Generator of canonical divisor}. On $U_{1,1}$, we use the representation of $\omega$ on $D(\Delta_{\{1,...,k \vert 1,...,k\}})$ and we have
\begin{flalign*}
	\pi_{0}^{*}(\bigwedge_{x_{i,j} \in \mathcal{S}_{\{1,...,k\vert 1,...,k\}}}dx_{i,j})=&\theta^{k(2r-k)-1}d\theta \wedge dy_{1,2} \wedge ... \wedge dy_{1,r} \\
	&\wedge dy_{2,1} \wedge ... \wedge dy_{r,1} \wedge (\bigwedge_{f_{i,j} \in \mathcal{S}_{\{2,...,k \vert 2,...,k\}}}df_{i,j});
\end{flalign*}
\begin{flalign*}
	\pi_{0}^{*}\Delta_{\{1,...,k\vert 1,...,k\}}=\theta^{k}\cdot det\begin{pmatrix}
		f_{2,2} & \dots & f_{2,r} \\
		\vdots  & \ddots & \vdots \\
		f_{r,2} & \dots & f_{r,r}
	\end{pmatrix}.
\end{flalign*}
On other charts the calculation is similar, so we have the relative canonical divisor of $\pi_{0}|_{\tilde{D^{k}}}$ is $(kr-1)\tilde{E_{0}}$. In subsequent  blowups, we can simply replace $k$ and $r$ with $k-1$ and $r-1$, respectively,  to obtain the desired result. In conclusion, we have
\begin{flalign*}
	K_{\tilde{\pi}}=\sum_{i=0}^{k-1}((k-i)(r-i)-1)\tilde{E_{i}}.
\end{flalign*}  
\end{proof}

\section{Stringy E-functions of determinantal varieties}
In this section, we derive the explicit formula for the stringy E-functions of determinantal varieties. Our primary approach involves expressing it through a form of integration utilizing the $\mathbb{Q}$-Gorenstein measure, as outlined in the following proposition.
\begin{proposition}\label{formula of stringy E-function via integration}
Suppose $X$ is a $\mathbb{Q}$-Gorenstein complex variety of pure dimension $d$ with at worst log terminal singularity, then the stringy E-function of $X$ is given by
\begin{flalign*}
	E_{st}(X)=E(\mu^{\mathbb{Q}-Gor}(J_{\infty}(X))).
\end{flalign*} 
\end{proposition}
\begin{proof}
Note that  this proposition can be derived  from the proof of Theorem $4.1$ in \cite{Stringy_E-function_via_Embdedded_Q-resolution},  however,  we  provide a proof here for  completeness.

Take a log resolution of $X$, $\phi: Y \longrightarrow X$. By Definition \ref{def_of_stringy_E-function}, we can rewrite the definition of stringy E-function into the integration as follows
\begin{flalign*}
	E_{st}(X)=E\big( \int_{J_{\infty}(Y)}\mathbb{L}^{-\mathrm{ord}_{t}K_{\phi}}\,d\mu_{J_{\infty}(Y)}\big).
\end{flalign*}
Since $Y$ is smooth, by Remark \ref{two_measure_coincide}, the $\mathbb{Q}$-Gorenstein measure and the usual measure on $J_{\infty}(Y)$ coincide, so $$E_{st}(X)=E( \int_{J_{\infty}(Y)}\mathbb{L}^{-\mathrm{ord}_{t}K_{\phi}}\,d\mu_{J_{\infty}(Y)}^{\mathbb{Q}\mathrm{-Gor}}).$$ Using the change of variables formula in \cite{Motivic_Zeta_Function_on_Q_Gorenstein_Varieties_Leon_Martin_Veys_Viu_Sos}(Theorem $2$), we get
\begin{flalign*}
	E_{st}(X)&=E(\int_{J_{\infty}(Y)}\mathbb{L}^{-\mathrm{ord}_{t}K_{\phi}}\,d\mu_{J_{\infty}(Y)}^{\mathbb{Q}\mathrm{-Gor}}) = E(\int_{J_{\infty}(X)}1\,d\mu_{J_{\infty}(X)}^{\mathbb{Q}\mathrm{-Gor}}).
\end{flalign*}
\end{proof}

From now on we assume $X=D^{k} \subset \mathcal{M}=\mathbb{C}^{r^{2}}$ is the determinantal variety defined by $(k+1)$ minors. By the definition of $\mu^{\mathbb{Q}-Gor}$, we need to calculate $\mathrm{ord}_{t}\omega_{X}(\psi)$ for every $\psi \in J_{\infty}(X)$. We fix an arc $\psi \in J_{\infty}(X)$, and we write $\psi=(x_{i,j})_{r \times r}$ as a matrix of arcs, where $x_{i,j}=x_{i,j}^{(0)}+x_{i,j}^{(1)}t+...$ is formal power series. By Proposition \ref{Generator of canonical divisor}, $\mathrm{ord}_{t}\omega_{X}(\psi)$ is given by $-(r-k)\mathrm{ord}_{t}D^{k-1}(\psi)$, where $\mathrm{ord}_{t}D^{k-1}(\psi)$ denotes the order of the ideal of $k$-minors under $\psi$. Note that this order is invariant under the action of $J_{\infty}(G)$, so it suffices to consider the orbit of $\psi$. By Proposition \ref{Orbit decomposition}, there is $\lambda=(\lambda_{1},...,\lambda_{r})$ such that $\infty=\lambda_{1}=...=\lambda_{r-k} \ge \lambda_{r-k+1} \ge ... \ge \lambda_{r}$ and $\psi \in \mathcal{C}_{\lambda}$, so we have $\mathrm{ord}_{t}\omega_{X}(\psi)=-(r-k)(\lambda_{r-k+1}+...+\lambda_{r})$. Moreover, for every $\psi \in \mathcal{C}_{\lambda}$, the order is the same, so we have
\begin{flalign*}
	E_{st}(X)&=E(\mu^{\mathbb{Q}-Gor}(J_{\infty}(X))) \\
             &=E(\int_{J_{\infty}(X)}\mathbb{L}^{-\mathrm{ord}_{t}\omega_{X}}\mathrm{d}\mu_{J_{\infty}(X)})  \\
             &=E(\sum_{\lambda_{r-k+1} \ge ... \ge \lambda_{r}}\mathbb{L}^{(r-k)(\lambda_{r-k+1}+...+\lambda_{r})}\mu_{J_{\infty}(X)}(\mathcal{C}_{\lambda})).	
\end{flalign*}

Note that if $\lambda_{r-k+1}=\infty$, we have $\mathcal{C}_{\lambda} \subset J_{\infty}(D^{k-1})$. Since $D^{k-1} \subset D^{k}$ is a proper closed subset, it is well-known that $\mu_{J_{\infty}(X)}(J_{\infty}(D^{k-1}))=0$, so we only need to consider the case $\lambda_{r-k+1}<\infty$. 

Next, we compute $\mu_{J_{\infty}(X)}(\mathcal{C}_{\lambda})$. We will mainly follow the notation of \cite{Arcs_on_Determinantal_Variteties_Roi}, which has calculated $\mu_{J_{\infty}(\mathcal{M})}(\mathcal{C}_{\lambda})$ in Section $6$.

We take an integer $n$ sufficiently large such that $n>\lambda_{r-k+1}$, so $\delta_{\lambda}$ can be regarded as an element in $J_{n}(X)$. Let $\mathcal{C}_{\lambda,n}$ be $n$-th jet truncation of $\mathcal{C}_{\lambda}$, and it suffices to compute $\mathbb{L}^{-nk(2r-k)}[\mathcal{C}_{\lambda,n}]$. Similar to Theorem \ref{Orbit decomposition}, $\mathcal{C}_{\lambda,n}$ is an orbit of $J_{n}(X)$ under the action of $J_{n}(G)$. Let $H_{\lambda,n}$ be the stabilizer of $\delta_{\lambda}$ under this action, and we have $[\mathcal{C}_{\lambda,n}]=\frac{[J_{n}(G)]}{[H_{\lambda,n}]}$. It is clear that $[J_{n}(G)]=[\mathrm{GL}_{r}]^{2} \cdot \mathbb{L}^{r^{2}n}$, so it suffices to consider $[H_{\lambda,n}]$.

For any element $(g,h) \in H_{\lambda,n}$, we write $g=(g_{i,j})_{r \times r},h=(h_{i,j})_{r \times r}$ and $g_{i,j}=g_{i,j}^{(0)}+g_{i,j}^{(1)}t+...+g_{i,j}^{(n)}t^{n},h_{i,j}=h_{i,j}^{(0)}+h_{i,j}^{(1)}t+...+h_{i,j}^{(n)}t^{n}$, then we have
\begin{flalign*}
(g,h) \in H_{\lambda,n} \Leftrightarrow g \cdot \delta_{\lambda} \cdot h^{-1}&=\delta_{\lambda} \Leftrightarrow g \cdot \delta_{\lambda}=\delta_{\lambda} \cdot h   \\
\Leftrightarrow (g_{i,j}) \cdot \mathrm{diag}\{0,...,0,t^{\lambda_{r-k+1}},...,t^{\lambda_{r}}\}&=\mathrm{diag}\{0,...,0,t^{\lambda_{r-k+1}},...,t^{\lambda_{r}}\} \cdot (h_{i,j}).   
\end{flalign*}
Comparing the elements on the two sides, we  get
\begin{flalign*}
	g_{i,j} t^{\lambda_{j}}&=0, \ 1 \le i \le r-k, \ r-k+1 \le j \le r; \\
	h_{i,j} t^{\lambda_{i}}&=0, \ r-k+1 \le i \le r, \ 1 \le j \le r-k; \\
	g_{i,j} t^{\lambda_{j}}&=h_{i,j} t^{\lambda_{i}}, \ r-k+1 \le i,j \le r.
\end{flalign*}
Let $H \subset G$ be the truncation of $H_{\lambda,n}$ to the constant jet level, so $H$ can be defined from $G$ with the following equations
\begin{flalign*}
	g_{i,j}^{(0)}&=0, \ 1 \le i \le r-k, \ r-k+1 \le j \le r; \\
	h_{i,j}^{(0)}&=0, \ r-k+1 \le i \le r, \ 1 \le j \le r-k; \\
	g_{i,j}^{(0)}&=h_{i,j}^{(0)}, \ r-k+1 \le i,j \le r, \ \lambda_{i}=\lambda_{j};  \\
	g_{i,j}^{(0)}&=0, \ r-k+1 \le i,j \le r, \ \lambda_{i} > \lambda_{j}; \\
	h_{i,j}^{(0)}&=0, \ r-k+1 \le i,j \le r, \ \lambda_{i}<\lambda_{j}.
\end{flalign*} 
In order to express these condition, we need to recall some notions in the group theory, as in the following definitions:
\begin{definition}
Assume $0 < v_{1}<...<v_{j}<r$ are fixed integers. We call a chain of $\mathbb{C}$-vector spaces $V_{1} \subset V_{2} \subset ... \subset V_{j} \subset \mathbb{C}^{r}$ a flag if the dimension of $V_{i}$ is $v_{i}$, and $(v_{1},...,v_{j})$ is called the signature of the flag. Note that if the signature is fixed, then $\mathrm{GL}_{r}$ acts transitively on the set of flags, and we call the stabilizer of a flag parabolic subgroup of $\mathrm{GL}_{r}$. If $P$ is a parabolic subgroup, then $\mathrm{GL}_{r}/P$ parametrizes the flag of a given signature.
\end{definition}
\begin{definition}\label{Definition of P and L}
Suppose $\{e_{1},...,e_{r}\}$ is a basis for $\mathbb{C}^{r}$, and $\lambda=(d_{1},...,d_{1},d_{2},...,d_{2},...,d_{j+1},...,d_{j+1})$ is a vector with $a_{i}$ many $d_{i}$, $d_{1}>...>d_{j+1} \ge 0$ and $a_{1}+...+a_{j+1}=r$. We set $v_{i}=a_{1}+...+a_{i}$ for $1 \le i \le j+1$, and $V_{i}=\mathrm{Span}\{e_{1},...,e_{v_{i}}\}$, then $V_{1} \subset... \subset V_{j} \subset V_{j+1}=\mathbb{C}^{r}$ is a flag. We denote the stabilizer of this flag $P_{\lambda}$ and we call it the parabolic subgroup associated to $\lambda$. $L_{\lambda}:=\mathrm{GL}_{a_{1}} \times \mathrm{GL}_{a_{2}} \times ... \times \mathrm{GL}_{a_{j+1}}$ is a subgroup of $P_{\lambda}$, and we call it the Levi factor associated to $\lambda$.
\end{definition}

If we fix a basis $\{e_{1},...,e_{r}\}$ as in the definition above and express the elements in $\mathrm{GL}_{r}$, $P_{\lambda}$ in matrices with respect to this basis, $P_{\lambda}$ will be those upper triangular matrices in blocks. The following example describes this explicitly.
\begin{example}
Suppose $r=5$, $\lambda=(3,3,1,1,0)$, $a_{1}=2,a_{2}=2,a_{3}=1$, then $v_{1}=2,v_{2}=4,v_{3}=5$. $P_{\lambda}$ and $L_{\lambda}$ are consisted of the matrices of the form
\begin{flalign*}
	P_{\lambda}=\begin{pmatrix}
		* & * & * & * & * \\
		* & * & * & * & * \\
		0 & 0 & * & * & * \\
		0 & 0 & * & * & * \\
		0 & 0 & 0 & 0 & * 
	\end{pmatrix}, L_{\lambda}=\begin{pmatrix}
	* & * & 0 & 0 & 0 \\
	* & * & 0 & 0 & 0 \\
	0 & 0 & * & * & 0 \\
	0 & 0 & * & * & 0 \\
	0 & 0 & 0 & 0 & * 
	\end{pmatrix}.
\end{flalign*}
\end{example}

Now we come back to the equations of $H$. The element of $H$ consists of two matrices $(g_{i,j}^{(0)}),(h_{i,j}^{(0)})$. $(g_{i,j}^{(0)})$ is lower triangular in blocks while $(h_{i,j}^{(0)})$ is upper triangular in blocks, and their blocks in the diagonal coincide except for the first one, as in the following picture shows.
\begin{flalign*}
	(g_{i,j}^{(0)})=\begin{pmatrix}
		A_{1} & 0 & 0 & \cdots & 0 \\
		* & B_{1} & 0 & \cdots & 0 \\
		* & * & B_{2} & \ddots & \vdots \\
		\vdots & \vdots & \vdots & \ddots & 0 \\
		* & * & * & * & B_{s} 
	\end{pmatrix}, (h_{i,j}^{(0)})=\begin{pmatrix}
		A_{2} & * & * & \cdots & * \\
		0 & B_{1} & * & \cdots & * \\
		0 & 0 & B_{2} & \cdots & * \\
		\vdots & \vdots & \ddots & \ddots & \vdots \\
		0 & 0 & \cdots & 0 & B_{s} 
	\end{pmatrix}.
\end{flalign*}
In the picture above, all the $*$ and $A_{1},A_{2},B_{1},...,B_{s}$ are matrices, and $A_{1},A_{2}$ are of the size $(r-k) \times (r-k)$. If we set $\lambda'=(n,...,n, \lambda_{r-k+1},...,\lambda_{r})$ and $\lambda''=(\lambda_{r-k+1},...,\lambda_{r})$, where $\lambda'$ has $r$ elements, then we have $H \cong (P_{\lambda'}^{T} \times P_{\lambda'})/L_{\lambda''}$. Here the action of $L_{\lambda''}$ is defined as 
\begin{flalign*}
P_{\lambda'}^{T} \times P_{\lambda'} \times L_{\lambda''} &\longrightarrow P_{\lambda'}^{T} \times P_{\lambda'}  \\
((g_{i,j}^{(0)}),(h_{i,j}^{(0)}),X) &\mapsto ((g_{i,j}^{(0)})X,(h_{i,j}^{(0)})X).	
\end{flalign*} 
The matrix $X$ is of the form
\begin{flalign*}
	\begin{pmatrix}
		Id & 0 & 0 & 0 & 0 \\
		0 & * & 0 & 0 & 0 \\
		0 & 0 & * & 0 & 0 \\
		0 & 0 & 0 & * & 0 \\
		0 & 0 & 0 & 0 & * 
	\end{pmatrix}
\end{flalign*}
and these $*$ are matrices with the same size as $B_{1},...,B_{s}$ so that $X$ can do the multiplication with $(g_{i,j}^{(0)})$ and $(h_{i,j}^{(0)})$.

Suppose we have determined the value of $g_{i,j}^{(0)}$ and $h_{i,j}^{(0)}$, then we consider the fiber over points of $H$. Noticed that the equations involving those $g_{i,j}^{(l)},h_{i,j}^{(l)}$ are linear equations, and every variable only appears one time in these equation, so the fiber will be a linear subspace of $\mathbb{C}^{r^{2}n}$, and the co-dimension will be the number of equations, which is
\begin{flalign*}
	&\sum_{i=r-k+1}^{r}(n-\lambda_{i})(r-k)+\sum_{j=r-k+1}^{r}(n-\lambda_{j})(r-k)+\sum_{r-k+1 \le i,j \le r}(n-\mathrm{min}\{\lambda_{i},\lambda_{j}\}) \\
	=&nk(2r-k)-\sum_{i=r-k+1}^{r}(2i-1)\lambda_{i}.
\end{flalign*}

In conclusion, we have
\begin{flalign*}
	\mu_{J_{\infty}(X)}(\mathcal{C}_{\lambda})=&\mathbb{L}^{-n(2r-k)k}\frac{[\mathrm{GL}_{r}]^{2} \cdot \mathbb{L}^{r^{2}n}}{[H_{\lambda,n}]} \\
	=&\mathbb{L}^{-n(2r-k)k}\frac{[\mathrm{GL}_{r}]^{2} \cdot \mathbb{L}^{r^{2}n} \cdot [L_{\lambda''}]}{\mathbb{L}^{r^{2}n-nk(2r-k)+\sum_{i=r-k+1}^{r}\lambda_{i}} \cdot [P_{\lambda'}]^{2}}  \\
	=&[\mathrm{GL}_{r}/P_{\lambda'}]^{2} \cdot [L_{\lambda''}] \cdot \mathbb{L}^{-\sum_{i=r-k+1}^{r}(2i-1)\lambda_{i}},
\end{flalign*}
and
\begin{flalign*}
	E_{st}(X)=E(\sum_{\lambda_{r-k+1} \ge ... \ge \lambda_{r}}[\mathrm{GL}_{r}/P_{\lambda'}]^{2} \cdot [L_{\lambda''}] \cdot \mathbb{L}^{\sum_{i=r-k+1}^{r}(r-k-2i+1)\lambda_{i}}).
\end{flalign*}
Noticed that if we write $\lambda'$ as the form in Definition \ref{Definition of P and L}, $[\mathrm{GL}_{r}/P_{\lambda'}]^{2} \cdot [L_{\lambda''}]$ only depends on the $a_{i}$ rather than $d_{i}$, so we want to sum all those $\lambda'$ having the same $a_{i}$. Now we set
\begin{flalign*}
	\lambda_{r}&=b_{r};  \\
	\lambda_{r-1}&=b_{r}+b_{r-1}; \\
	&\ \ \vdots  \\
	\lambda_{r-k+1}&=b_{r}+...+b_{r-k+1}.
\end{flalign*}
It is clear that vector $(\lambda_{r-k+1},...,\lambda_{r})$ and $(b_{r-k+1},...,b_{r})$ are determined by each other. For any $I \subset \{r-k+1,...,r-1\}$, let $\Omega_{I}:=\{(b_{r-k+1},...,b_{r})\;\vert\; \forall i \in I, b_{i}=0, \forall i \notin I, b_{i} \neq 0\}$, then for any two different $\lambda$, if their corresponding $(b_{r-k+1},...,b_{r})$ belong to the same $\Omega_{I}$, then the value of $[\mathrm{GL}_{r}/P_{\lambda'}]^{2} \cdot [L_{\lambda''}]$ are the same. Let $I_{r}^{c}=\{r-k+1,...,r\}-I$, and assume $I_{r}^{c}=\{i_{1},...,i_{l}\}$, where $i_{1}<...<i_{l}$, then 
\begin{flalign*}
	[L_{\lambda}'']&=\prod_{j=1}^{l}[\mathrm{GL}_{i_{j}-i_{j-1}}]; \\
	[\mathrm{GL}_{r}/P_{\lambda'}]&= \prod_{j=1}^{l}[G(i_{j}-i_{j-1},i_{j})],
\end{flalign*}
where $i_{0}=r-k$  and $G(u,v)$ denotes the Grassmannian of $u$-dimensional subspace in the $v$-dimensional vector space. 

If we express $\mathbb{L}^{\sum_{i=r-k+1}^{r}(r-k-2i+1)\lambda_{i}}$ with the information of $b_{r-k+1},...,b_{r}$, we have
\begin{flalign*}
	\mathbb{L}^{\sum_{i=r-k+1}^{r}(r-k-2i+1)\lambda_{i}}  
	=\mathbb{L}^{\sum_{i=r-k+1}^{r}-i(i-r+k)b_{i}}.
\end{flalign*}

We define $d(I,i_{j}):=i_{j}-i_{j-1}$ for $i_{j} \in I_{r}^{c}$, and we can get
\begin{flalign*}
	E_{st}(X)=E(\sum_{I \subset \{r-k+1,...,r-1\}}\sum_{b(\lambda) \in \Omega_{I}}\prod_{i \in I_{r}^{c}}[\mathrm{GL}_{d(I,i)}][G(d(I,i),i)]^{2}\mathbb{L}^{\sum_{i=r-k+1}^{r}-i(i-r+k)b_{i}}).
\end{flalign*}

Note that
\begin{flalign*}
	&\sum_{b(\lambda) \in \Omega_{I}}\mathbb{L}^{\sum_{i=r-k+1}^{r}-i(i-r+k)b_{i}}  \\
	=&\sum_{b(\lambda) \in \Omega_{I}}\mathbb{L}^{\sum_{i \in I_{r}^{c}}-i(i-r+k)b_{i}}  \\
	=&(\prod_{i \in \{r-k+1,...,r-1\}-I}\sum_{b_{i}=1}^{\infty}\mathbb{L}^{-b_{i}i(i-r+k)})\sum_{b_{r}=0}^{\infty}\mathbb{L}^{-b_{r}kr}  \\
	=&(\prod_{i \in \{r-k+1,...,r-1\}-I}\frac{1}{\mathbb{L}^{i(i-r+k)}-1}) \cdot \frac{\mathbb{L}^{kr}}{\mathbb{L}^{kr}-1}  \\
	=&\mathbb{L}^{kr}\prod_{i \in I_{r}^{c}}\frac{1}{\mathbb{L}^{i(i-r+k)}-1}.
\end{flalign*}

In conclusion, we have the following Proposition.
\begin{proposition}\label{complicated formula}
The stringy E-function of $D^{k}$ is given by
\begin{flalign*}
	E_{st}(D^{k})=E(\mathbb{L}^{kr}\sum_{I \subset \{r-k+1,...,r-1\}}\prod_{i \in I_{r}^{c}}\frac{1}{\mathbb{L}^{i(i-r+k)}-1}[\mathrm{GL}_{d(I,i)}][G(d(I,i),i)]^{2}).
\end{flalign*}
\end{proposition}

\begin{remark}
When $k=r-1$, we have $$E_{st}(D^{r-1})=E(\mathbb{L}^{r(r-1)}\sum_{I \subset \{2,...,r-1\}}\prod_{i \in I_{r}^{c}}\frac{1}{\mathbb{L}^{i(i-1)}-1}[\mathrm{GL}_{d(I,i)}][G(d(I,i),i)]^{2}).$$ On the other hand, since $D^{r-1}$ is defined by a single equation, one can use Proposition $2.3$ in \cite{Stringy_E-function_of_Hypersurface} and the motivic zeta function of this equation to compute it. The latter has been given in Section $6$ in \cite{Arcs_on_Determinantal_Variteties_Roi}, which is
\begin{flalign*}
	Z_{D^{r-1}}(T)=\mathbb{L}^{r^{2}}T^{-r}\sum_{I \subset \{1,...,r-1\}}\prod_{i \in I_{r}^{c}}\frac{1}{\mathbb{L}^{i^{2}}T^{-i}-1}[\mathrm{GL}_{d(I,i)}][G(d(I,i),i)]^{2}.
\end{flalign*}
Note that in the notation of \cite{Arcs_on_Determinantal_Variteties_Roi}, $i_{0}=0$ in the calculation of $d(I,i)$. After taking the hodge E-function and residue at $uv$ as in Proposition $2.3$ in \cite{Stringy_E-function_of_Hypersurface}, one can verify that the stringy E-function obtained by these two methods coincides.
\end{remark}

\begin{remark}\label{k=1}
If we take $k=1$, we can get $E_{st}(D^{1})=(uv)^{r}\frac{(uv)^{r}-1}{uv-1}$. On the other hand, by Proposition \ref{resolution of Dk}, the resolution of $D^{1}$ can be obtained by one blowup, so we can use the definition of stringy E-function to compute it, which is
\begin{flalign*}
	E_{st}(D^{1})=&E(\tilde{D^{1}}-\tilde{E_{0}})+E(\tilde{E_{0}}) \cdot \frac{uv-1}{(uv)^{r}-1} \\
	=&E(D^{1})-1+E(\tilde{E_{0}})\cdot \frac{uv-1}{(uv)^{r}-1}.
\end{flalign*}
One can check that $E(D^{1})=1+\frac{((uv)^{r}-1)^{2}}{uv-1}$ and $E(\tilde{E_{0}})=\frac{E(D^{1})-1}{uv-1}$, and we can get the same result through this way.
\end{remark}

Finally, we will prove that this expression is exactly the E-function of a Grassmannian. Before that, we need some preparation about some special varieties in the Grothendieck ring.
\begin{lemma}[Example $2.4.4$ and $2.4.5$ in \cite{Motivic_Integration_book}]\label{Grassmannian}
For integer $0<d \le k$, we have:

(i)The general linear group $[\mathrm{GL}_{d}]=\mathbb{L}^{d(d-1)/2}(\mathbb{L}^{d}-1)(\mathbb{L}^{d-1}-1) \cdots (\mathbb{L}-1)$.

(ii)The Grassmannian $[G(d,k)]=\prod_{j=1}^{d}\frac{\mathbb{L}^{j+k-d}-1}{\mathbb{L}^{j}-1}=\sum_{0 \le \lambda_{1} \le ... \le \lambda_{d} \le k-d}\mathbb{L}^{\lambda_{1}+...+\lambda_{d}}$.

(iii)Suppose $V=\mathbb{C}^{k}$ and $U(d,k)$ is the open subscheme of $V^{d}$ consisting of linearly independent vectors $(v_{1},...,v_{d})$ with $v_{i} \in V$, then $[U(d,k)]=(\mathbb{L}^{k}-\mathbb{L}^{d-1})(\mathbb{L}^{k}-\mathbb{L}^{d-2}) \cdots (\mathbb{L}^{k}-1)$.
\end{lemma}

\begin{theorem}\label{formula}
\begin{flalign*}
	E_{st}(D^{k})=E(\mathbb{L}^{kr}\prod_{i=1}^{k}\frac{\mathbb{L}^{i+r-k}-1}{\mathbb{L}^{i}-1})=E(\mathbb{L}^{kr} \cdot G(k,r)).
\end{flalign*}
\end{theorem}
\begin{proof}
It suffices to prove the following equation
\begin{flalign*}
	A(k,r):=\sum_{I \subset \{r-k+1,...,r-1\}}\prod_{i \in I_{r}^{c}}\frac{1}{\mathbb{L}^{i(i-r+k)}-1}[\mathrm{GL}_{d(I,i)}][G(d(I,i),i)]^{2}=\prod_{i=1}^{k}\frac{\mathbb{L}^{i+r-k}-1}{\mathbb{L}^{i}-1}.
\end{flalign*}
We use induction on $k+r$. When $k=1$, this equation holds due to Remark \ref{k=1}. On the other hand, by Lemma \ref{Grassmannian}(i)(ii), we have $[\mathrm{GL}_{d}]=\mathbb{L}^{d(d-1)/2}(\mathbb{L}^{d}-1)(\mathbb{L}^{d-1}-1) \cdots (\mathbb{L}-1)$ and $[G(d,i)]=\prod_{j=1}^{d}\frac{\mathbb{L}^{j+i-d}-1}{\mathbb{L}^{j}-1}$, so  $$A(k,r)=\sum_{i_{0}=r-k<i_{1}<...<i_{l-1}<i_{l}=r}\prod_{j=1}^{l}\frac{1}{\mathbb{L}^{i_{j}(i_{j}-r+k)}-1} \frac{(\mathbb{L}^{i_{j-1}+1}-1)^{2} \cdots (\mathbb{L}^{i_{j}}-1)^{2}}{(\mathbb{L}-1) \cdots (\mathbb{L}^{i_{j}-i_{j-1}}-1)} \mathbb{L}^{\frac{(i_{j}-i_{j-1})(i_{j}-i_{j-1}-1)}{2}}.$$ When examining the summation, we can treat $i_{l-1}$ as a new variable. By substituting  $r$ and $k$ with  $i_{l}$ and $k-r+i_{l}$, respectively, we obtain the following relation
\begin{flalign*}
A(k,r)=&\sum_{m=r-k+1}^{r-1}A(k+m-r,m)\frac{1}{\mathbb{L}^{kr}-1} \frac{(\mathbb{L}^{m+1}-1)^{2} \cdots (\mathbb{L}^{r}-1)^{2}}{(\mathbb{L}-1) \cdots (\mathbb{L}^{r-m}-1)} \mathbb{L}^{\frac{(r-m)(r-m-1)}{2}} \\
&+\frac{1}{\mathbb{L}^{kr}-1} \frac{(\mathbb{L}^{r-k+1}-1)^{2} \cdots (\mathbb{L}^{r}-1)^{2}}{(\mathbb{L}-1) \cdots (\mathbb{L}^{k}-1)} \mathbb{L}^{\frac{k(k-1)}{2}}.
\end{flalign*} 
Using induction hypothesis, it suffices to prove
\begin{flalign*}
&\sum_{m=r-k+1}^{r-1}\frac{(\mathbb{L}^{r-k+1}-1) \cdots (\mathbb{L}^{m}-1)}{(\mathbb{L}-1) \cdots (\mathbb{L}^{k+m-r}-1)}\frac{1}{\mathbb{L}^{kr}-1} \frac{(\mathbb{L}^{m+1}-1)^{2} \cdots (\mathbb{L}^{r}-1)^{2}}{(\mathbb{L}-1) \cdots (\mathbb{L}^{r-m}-1)} \mathbb{L}^{\frac{(r-m)(r-m-1)}{2}} \\
&+\frac{1}{\mathbb{L}^{kr}-1} \frac{(\mathbb{L}^{r-k+1}-1)^{2} \cdots (\mathbb{L}^{r}-1)^{2}}{(\mathbb{L}-1) \cdots (\mathbb{L}^{k}-1)} \mathbb{L}^{\frac{k(k-1)}{2}}=\frac{(\mathbb{L}^{r-k+1}-1) \cdots (\mathbb{L}^{r}-1)}{(\mathbb{L}-1) \cdots (\mathbb{L}^{k}-1)}.
\end{flalign*}
After cancellation, we only need to prove
\begin{flalign*}
	\mathbb{L}^{kr}-1=\sum_{m=r-k}^{r-1}(\mathbb{L}^{k}-\mathbb{L}^{r-m-1})(\mathbb{L}^{k}-\mathbb{L}^{r-m-2}) \cdots (\mathbb{L}^{k}-1) \cdot \frac{(\mathbb{L}^{m+1}-1) \cdots (\mathbb{L}^{r}-1)}{(\mathbb{L}-1) \cdots (\mathbb{L}^{r-m}-1)}.
\end{flalign*}
From Lemma \ref{Grassmannian}(ii)(iii),  we obtain the following identities $$\frac{(\mathbb{L}^{m+1}-1) \cdots (\mathbb{L}^{r}-1)}{(\mathbb{L}-1) \cdots (\mathbb{L}^{r-m}-1)}=[G(m,r)]$$ and  $$(\mathbb{L}^{k}-\mathbb{L}^{r-m-1})(\mathbb{L}^{k}-\mathbb{L}^{r-m-2}) \cdots (\mathbb{L}^{k}-1)=[U(r-m,k)],$$ where $U(r-m,k)$ denotes the open subscheme of $V^{r-m}$ consisting of linearly independent vectors $(v_{1},...,v_{r-m})$ with  $v_{i} \in V=\mathbb{C}^{k}$. Let $M(k,r)_{r-m}$ be the space of linear maps from $\mathbb{C}^{r}$ to $\mathbb{C}^{k}$ with rank $r-m$, where $r-k \le m \le r-1$. For any $\phi \in M(k,r)_{r-m}$, consider the map
\begin{flalign*}
	\Phi:M(k,r)_{r-m} &\longrightarrow G(m,r) \\
	        \phi  &\mapsto \mathrm{Ker}(\phi).
\end{flalign*}
If we fix $\mathrm{Ker}(\phi)$ and we take a basis of the complementary space of it in $\mathbb{C}^{r}$, then $\phi$ is determined by the image of this basis, which is a set of linearly independent vectors $(v_{1},...,v_{r-m})$. This tells us the fibers of $\Phi$ are $U(r-m,k)$, so we have $[M(k,r)_{r-m}]=[G(m,r)] \cdot [U(r-m,k)]$. Note that the space of linear map with rank $0$ is a point and Grothendieck group of the whole space of linear maps from $\mathbb{C}^{r}$ to $\mathbb{C}^{k}$ is $\mathbb{L}^{kr}$, so we get the following equation
\begin{flalign*}
	\mathbb{L}^{kr}=1+\sum_{m=r-k}^{r-1}[G(m,r)][U(r-m,k)].
\end{flalign*}
This gives the required statement.  
\end{proof}
As direct applications, we have the following corollary.
\begin{corollary}
The stringy Euler number of $D^{k}$, defined by $\mathrm{lim}_{u,v \rightarrow 1}E_{st}(D^{k})$, is $C^{k}_{r}$.
\end{corollary}

\section{Stringy E-functions of projective determinantal varieties}
In this section, we compute the stringy E-functions of the projectivization of determinantal varieties. To ensure clarity, we adopt a more precise notation. Let $\mathcal{M}_{r,s}=\mathbb{C}^{rs}$ be the space of $r \times s$ matrices with coordinate ring $\mathbb{C}[x_{i,j}]_{1 \le i \le r, 1 \le j \le s}$, and $\hat{\mathcal{M}}_{r,s}:=\mathbb{P}(\mathcal{M}_{r,s})=\mathbb{P}^{rs-1}$ is the projectivization of it. We assume $D^{k}_{r,s}$ is the subvariety of $\mathcal{M}_{r,s}$ defined by $(k+1)$-minors and $\hat{D}^{k}_{r,s} \subset \hat{\mathcal{M}}_{r,s}$ is projectivization of it. Suppose $U_{i,j}$ is the open set of $\hat{\mathcal{M}}_{r,s}$ such that $x_{i,j} \neq 0$, then these $U_{i,j}$ consists of an open covering of $\hat{\mathcal{M}}_{r,s}$. The following relation shows that locally $\hat{D}^{k}_{r,s}$ is still a determinantal variety.
\begin{proposition}\label{locally determinantal}
We have isomorphism $$U_{i,j}\cap \hat{D}^{k}_{r,s} \cong \mathbb{C}^{r-1} \times \mathbb{C}^{s-1} \times D^{k-1}_{r-1,s-1}.$$
\end{proposition}  
\begin{proof}
Note that this isomorphism can be obtained from Section $3.4$ in \cite{Multiplier_Ideals_of_Determinantal_Ideal}, but we will still give the construction of it for completeness. It suffices to prove the case when $i=j=1$. We use $y_{i,j}:=\frac{x_{i,j}}{x_{1,1}}$ as the affine coordinates in $U_{1,1}$ and we consider the following morphism.
\begin{flalign*}
	\phi: U_{i,j}\cap \hat{D}^{k}_{r,s} &\longrightarrow \mathbb{C}^{r-1} \times \mathbb{C}^{s-1} \times D^{k-1}_{r-1,s-1}  \\
	\begin{pmatrix}
		1       & y_{1,2} & \cdots & y_{1,s} \\
		y_{2,1} & y_{2,2} & \cdots & y_{2,s} \\
		\vdots  & \vdots  & \ddots & \vdots  \\
		y_{r,1} & y_{r,2} & \cdots & y_{r,s}
	\end{pmatrix} &\mapsto (y_{2,1},...,y_{r,1}) \times (y_{1,2},...,y_{1,s}) \times (z_{i,j}),
\end{flalign*} 
where $(z_{i,j})$ denotes the $(r-1) \times (s-1)$ matrix with entries $z_{i,j}=y_{i,j}-y_{i,1}y_{1,j}(2 \le i \le r, 2 \le j \le s)$. This morphism $\phi$ actually cancels these $y_{i,1}$ and $y_{1,j}$ by row and column transformations, so the rank of $(y_{i,j})$ is less or equal to $k$ is equivalent to $r((z_{i,j})) \le k-1$. One can easily give the inverse of $\phi$ and verify that this is an isomorphism.
\end{proof}

Note that if $k=1$, we have $\hat{D}^{k}_{r,s} \cong \mathbb{P}^{r-1} \times \mathbb{P}^{s-1}$, so we will focus on $k \ge 2$. Combining with the result in the last section, we know that $\hat{D}^{k}_{r,s}$ is Gorenstein if and only if $r=s$, and in this case it has canonical singularity, so we will assume $r=s$ in our later study. 

We will employ the same method as in the previous section to compute the stringy E-function. Prior to this, we need to extend Proposition \ref{Orbit decomposition} to the projective case. Firstly we characterize $J_{\infty}(\hat{\mathcal{M}}_{r,r})$ and $J_{\infty}(\hat{D}^{k}_{r,r})$. $J_{\infty}(\hat{\mathcal{M}}_{r,r})=J_{\infty}(\mathbb{P}^{r^{2}-1})$ can be regarded as
\begin{flalign*}
	\lbrace \begin{pmatrix}
		x_{1,1}^{(0)}+x_{1,1}^{(1)}t+... & \cdots & x_{1,r}^{(0)}+x_{1,r}^{(1)}t+...   \\
		\vdots & \ddots & \vdots \\
		x_{r,1}^{(0)}+x_{r,1}^{(1)}t+... & \cdots & x_{r,r}^{(0)}+x_{r,r}^{(1)}t+...
	\end{pmatrix}\; \vert \; \exists \ i,j , \ x_{i,j}^{(0)} \neq 0 \rbrace / \thicksim,
\end{flalign*}
where the equivalence $\thicksim$ is given by multiplying an invertible power series $f$ at all entries
\begin{flalign*}
	\begin{pmatrix}
		x_{1,1}^{(0)}+x_{1,1}^{(1)}t+... & \cdots & x_{1,r}^{(0)}+x_{1,r}^{(1)}t+...   \\
		\vdots & \ddots & \vdots \\
		x_{r,1}^{(0)}+x_{r,1}^{(1)}t+... & \cdots & x_{r,r}^{(0)}+x_{r,r}^{(1)}t+...
	\end{pmatrix} \thicksim \begin{pmatrix}
	f\cdot (x_{1,1}^{(0)}+x_{1,1}^{(1)}t+...) & \cdots & f \cdot (x_{1,r}^{(0)}+x_{1,r}^{(1)}t+...)   \\
	\vdots & \ddots & \vdots \\
	f \cdot (x_{r,1}^{(0)}+x_{r,1}^{(1)}t+...) & \cdots & f \cdot (x_{r,r}^{(0)}+x_{r,r}^{(1)}t+...)
	\end{pmatrix}.
\end{flalign*}
By Proposition \ref{Orbit decomposition}, we have the orbit decomposition $$J_{\infty}(\mathcal{M}_{r,r})=\bigcup_{\lambda}\mathcal{C}_{\lambda}$$ with respect to $J_{\infty}(\mathrm{GL}_{r} \times \mathrm{GL}_{r})$, where $\lambda=(\lambda_{1},...,\lambda_{r})$ with $\infty \ge \lambda_{1} \ge ... \ge \lambda_{r} \ge 0$. Note that the action of this group is compatible with the equivalence above, so we have similar decomposition 
\begin{flalign*}
J_{\infty}(\hat{\mathcal{M}}_{r,r})=\bigcup_{ \lambda_{r}=0}\mathcal{C}_{\lambda} / \thicksim.
\end{flalign*}
Again by Proposition \ref{Orbit decomposition}, $J_{\infty}(D^{k}_{r,r})$ can be written as
\begin{flalign*}
	J_{\infty}(D^{k}_{r,r})=\bigcup_{\lambda_{1}=...=\lambda_{r-k}=\infty}\mathcal{C}_{\lambda},
\end{flalign*}
and we have
\begin{flalign*}
	J_{\infty}(\hat{D}^{k}_{r,r})=\bigcup_{\lambda_{1}=...=\lambda_{r-k}=\infty,\lambda_{r}=0}\mathcal{C}_{\lambda}/ \thicksim.
\end{flalign*}
Now we start to calculate the stringy E-function of $\hat{D}^{k}_{r,r}$. By applying Proposition \ref{formula of stringy E-function via integration}, we obtain
\begin{flalign*}
	E_{st}(\hat{D}^{k}_{r,r})=E(\int_{J_{\infty}(\hat{D}^{k}_{r,r})}\mathbb{L}^{-\mathrm{ord}_{t}\omega_{\hat{D}^{k}_{r,r}}}d\mu_{J_{\infty}(\hat{D}^{k}_{r,r})}).
\end{flalign*}
By Proposition \ref{locally determinantal} and Proposition \ref{Generator of canonical divisor}, we can characterize $\omega_{\hat{D}^{k}_{r,r}}$ in a similar way. So for an arc $\psi \in \mathcal{C}_{\lambda}/\thicksim$ with $\lambda_{1}=...=\lambda_{r-k}=\infty$ and $\lambda_{r}=0$, we have $-\mathrm{ord}_{t}\omega_{\hat{D}^{k}_{r,r}}(\psi)=(r-k)(\lambda_{r-k+1}+...+\lambda_{r-1})$. This tells us
\begin{flalign*}
	E_{st}(\hat{D}^{k}_{r,r})=E(\sum_{\lambda_{1}=...=\lambda_{r-k}=\infty,\lambda_{r}=0}\mu_{J_{\infty}(\hat{D}^{k}_{r,r})}(\mathcal{C}_{\lambda}/\thicksim) \cdot  \mathbb{L}^{(r-k)(\lambda_{r-k+1}+...+\lambda_{r-1})}) .
\end{flalign*}
Similar to the proof of Proposition \ref{complicated formula}, we choose a number $n$ sufficiently large so that we can consider $\mathcal{C}_{\lambda}$ in the $n$-th jet level, and we have
\begin{flalign*}
	\mu_{J_{\infty}(\hat{D}^{k}_{r,r})}(\mathcal{C}_{\lambda}/\thicksim)=\mathbb{L}^{-n(k(2r-k)-1)}[\mathcal{C}_{\lambda,n}/\thicksim].
\end{flalign*}
Note that $[\mathcal{C}_{\lambda,n}/\thicksim]=\frac{[\mathcal{C}_{\lambda,n}]}{(\mathbb{L}-1)\mathbb{L}^{n}}$, so we have $\mu_{J_{\infty}(\hat{D}^{k}_{r,r})}(\mathcal{C}_{\lambda}/\thicksim)=\frac{1}{\mathbb{L}-1}\cdot \mu_{J_{\infty}({D}^{k}_{r,r})}(\mathcal{C}_{\lambda})$, which implies
\begin{flalign*}
	E_{st}(\hat{D}^{k}_{r,r})=E(\frac{1}{\mathbb{L}-1}\sum_{\lambda_{1}=...=\lambda_{r-k}=\infty,\lambda_{r}=0}\mu_{J_{\infty}(D^{k}_{r,r})}(\mathcal{C}_{\lambda}) \cdot  \mathbb{L}^{(r-k)(\lambda_{r-k+1}+...+\lambda_{r-1}+\lambda_{r})}).
\end{flalign*}
Again we use the proof and the notation of Proposition \ref{complicated formula}, to be precise, let $\lambda_{i}=b_{r}+b_{r-1}+...+b_{i}(r-k+1 \le i \le r)$, and for any set $I \subset \{r-k+1,...,r-1\}$, let $\Omega_{I}:=\{(b_{r-k+1},...,b_{r})\;\vert \;\forall i \in I, b_{i}=0, \forall i \notin I, b_{i} \neq 0\}$. Suppose $I^{c}_{r}:=\{r-k+1,...,r\}-I=\{i_{1},...,i_{l}\}$ such that $r-k=i_{0}<i_{1}<...<i_{l}$, and we set $d(I,i_{j}):=i_{j}-i_{j-1}$. Finally we will get
\begin{flalign*}
	E_{st}(\hat{D}^{k}_{r,r})=E(\sum_{I \subset \{r-k+1,...,r-1\}}\sum_{b(\lambda) \in \Omega_{I},b_{r}=0}\prod_{i \in I_{r}^{c}}[\mathrm{GL}_{d(I,i)}][G(d(I,i),i)]^{2}\mathbb{L}^{\sum_{i=r-k+1}^{r}-i(i-r+k)b_{i}}).
\end{flalign*}
Note that
\begin{flalign*}
	\sum_{b(\lambda) \in \Omega_{I},b_{r}=0}\mathbb{L}^{\sum_{i=r-k+1}^{r}-i(i-r+k)b_{i}}  
	=(\mathbb{L}^{kr}-1)\prod_{i \in I_{r}^{c}}\frac{1}{\mathbb{L}^{i(i-r+k)}-1}.
\end{flalign*}
Combining with the proof of Theorem \ref{formula}, we obtain the following result.
\begin{theorem}\label{formula of projective determinantal variety}
$\hat{D}^{k}_{r,r}$ has Gorenstein canonical singularity and the stringy E-function of it is
\begin{flalign*}
	E_{st}(\hat{D}^{k}_{r,r})=\frac{(uv)^{kr}-1}{(uv-1) \cdot (uv)^{kr}}\cdot E_{st}(D^{k}_{r,r})=E(\frac{\mathbb{L}^{kr}-1}{\mathbb{L}-1} \cdot \prod_{i=1}^{k}\frac{\mathbb{L}^{i+r-k}-1}{\mathbb{L}^{i}-1})=\frac{(uv)^{kr}-1}{uv-1} \cdot E([G(k,r)]).
\end{flalign*}
\end{theorem}

\begin{remark}
Although we suppose $k \ge 2$ in the discussion above, one can verify that this theorem still holds for $k=1$.
\end{remark}

\begin{corollary}\label{c5.4}
The non-negativity conjecture of stringy Hodge numbers holds for $\hat{D}^{k}_{r,r}$.
\end{corollary}
\begin{proof}
By Lemma \ref{Grassmannian}(ii), $[G(k,r)]$ can be written as a polynomial of $\mathbb{L}$ with positive coefficients, so the result follows from Theorem \ref{formula of projective determinantal variety}. 
\end{proof}

\begin{proof}[proof of Theorem \ref{th2}.]

Combining Theorem \ref{formula of projective determinantal variety} with Corollary \ref{c5.4}, we obtain Theorem \ref{th2}.
\end{proof}

As a corollary, we also derive the stringy Euler number of  $\hat{D}^{k}_{r,r}$.
\begin{corollary}
The stringy Euler number of $\hat{D}^{k}_{r,r}$ is $kr\cdot C^{k}_{r}$.
\end{corollary}

\end{document}